\documentclass[a4paper,12pt]{amsart}
\usepackage{kv1}

\begin{document}

\title{Physical measures for nonlinear random walks on interval}

\author{V.~Kleptsyn}
\address{Victor~Kleptsyn
\newline\hphantom{iii} CNRS, Institut de Recherche Mathematique de Rennes (IRMAR, UMR 6625 CNRS)}
\email{Victor.Kleptsyn@univ-rennes1.fr}
\thanks{V.~K. was supported in part by grants RFBR 10-01-00739, RFBR/CNRS 10-01-93115}

\author{{D.~Volk}}
\address{Denis~Volk
\newline\hphantom{iii} University of Rome ``Tor Vergata''
\newline\hphantom{iii} Institute for Information Transmission Problems, Russian Academy of Sciences}
\email{volk@mat.uniroma2.it}
\thanks{D.~V. was supported in part by grants RFBR 10-01-00739, RFBR/CNRS 10-01-93115, President's of Russia MK-2790.2011.1 and MK-7567.2013.1, PRIN, ``Young SISSA Scientists'', RFBR 12-01-31241-mol\_a, RFBR 13-01-00969-a}

\subjclass[2010]{Primary: 82B41, 82C41, 60G50. Secondary: 37C05, 37C20, 37C70, 37D45}

\keywords{Random walks, stationary measures, dynamical systems, attractors, partial hyperbolicity, skew products}

\begin{abstract}
A one-dimensional confined Nonlinear Random Walk is a tuple of $N$ diffeomorphisms of the unit interval driven by a probabilistic Markov chain. For generic such walks, we obtain a geometric characterization of their ergodic stationary measures and prove that all of them have negative Lyapunov exponents.

These measures appear to be probabilistic manifestations of physical measures for certain deterministic dynamical systems. These systems are step skew products over transitive subshifts of finite type (topological Markov chains) with the unit interval fiber.

For such skew products, we show there exist only finite collection of alternating attractors and repellers; we also give a sharp upper bound for their number. Each of them is a graph of a continuous map from the base to the fiber defined almost everywhere w.r.t. any ergodic Markov measure in the base. The orbits starting between the adjacent attractor and repeller tend to the attractor as $t \to +\infty$, and to the repeller as $t \to -\infty$. The attractors support ergodic hyperbolic physical measures.
\end{abstract}

\maketitle

\section{Introduction}  \label{sec:intro}


\subsection{Nonlinear Random Walks}

The classical one-dimensional discrete-time random walk on the real line is the process when every time you toss a coin and make a jump left or right at one unit depending on the toss outcome.

In this paper, we consider a far generalization of this construction. First, we allow more than two kinds of jumps, and let them depend on the previous jump. In other words, we let the jumps to be driven by a Markov chain with finitely many states. Second, we let the jumps to depend on the current position on the real line, and do this in a nonlinear way. From the point $x$ you make a jump to the point~$f_i(x)$, where $f_i$ are some
diffeomorphisms of the real line. Finally, we make a very important assumption which will play through all the paper: that the whole process is confined within a unit interval of the real line.

In these very general settings we study the probabilistic stationary measures of such processes. It turns out that for a generic process of this kind, there exist only finitely many ergodic stationary measures, and we give a sharp estimate on their number. Their geometry and relative positions on the unit interval can be described in terms of fixed points of the jump maps, see Theorem~\ref{t:random}. Moreover, their Lyapunov exponents are always negative, so the orbits converge on average for almost every random itinerary.

%


\subsection{Skew products}

Random walks have their deterministic twins, namely, the step skew product dynamical systems over topological Markov chains (subshifts of finite type), see Definition~\ref{d:step_and_mild}.
In our setting, these are skew products whose fibers are unit intervals, and fiber maps are $C^1$-diffeomorphisms onto the image. 

In the dynamical version of our main theorem we show that the dynamics of a generic step skew product with interval fibers can be described in relatively simple terms. This dynamics is in a sense similar to the cartesian product of the dynamics in the base and the dynamics of a single interval diffeomorphism. More precisely (see Theorem~\ref{t:markov_classification}), the phase space of such a skew product can be covered by finitely many absorbing and expelling strips. Any absorbing strip contains a unique attractor; any expelling strip contains a unique repeller. These attractors and repellers are \emph{bony graphs} of maps from the base to the fiber: each of them intersects almost every fiber (w.r.t. the Markov measure in the base) at a single point; other fibers are intersected at intervals. This feature is similar to \emph{porcupine horseshoes} discovered by {D\'{\i}az} and Gelfert in~\cite{Diaz2012}. In particular, we partially answer Question 1.6 from~\cite{Diaz2012} by showing that generically, the set of bones is meager. Almost every point of the phase space (w.r.t. the \emph{standard measure} which is the product of the Markov measure in the base and Lebesgue measure in the fiber, sf. Definition~\ref{d:std}) tends to one of the attractor graphs. When the time is reversed, almost every point either tends to one of the repeller graphs, or is eventually taken to a domain where the inverse map is not defined.

Finally, the closure of each attractor graph and of each repeller graph is the support of an ergodic invariant measure that projects to the Markov measure in the base. The attractor measures are physical (SRB); their basins of attraction contain subsets of full standard measure in the corresponding strips. The fiberwise Lyapunov exponents of the attracting and the repelling measures are strictly negative and strictly positive, respectively.

In~\cite{Kudryashov2010}, Kudryashov has recently constructed a robust example of a step skew product with interval fibers such that its attractor intersects some fibers at intervals (instead of points). The set of such fibers has continuum cardinality (in fact, it has a Cantor subset). Thus it is not possible to prove the stronger statement that the attractors and repellers intersect each fiber at a single point only.

\subsection{Background and motivations}

There are important reasons to study the associated skew product dynamical systems on their own, too. We provide a brief review below.

Among the dynamical systems, the open set of hyperbolic ones (also known as Axiom~A) is the best understood. For them, Smale's Spectral Decomposition theorem~\cite{Smale1967a} allows one to split the non-wandering set into a finite collection of locally maximal hyperbolic sets. These unit sets admit finite Markov partitions and, therefore, symbolic encodings. If any of them is an attractor, then it also carries an SRB (physical) measure~\cite{Sinai1972}, \cite{Ruelle1976}, \cite{Bowen2008}. Thus the hyperbolic systems have both good dynamical and statistical descriptions.

From the works of Abraham, Smale~\cite{Abraham1970} and Newhouse~\cite{Newhouse1974} we know that Axiom~A is not dense in $\dim \ge 2$. So the next step to understand generic systems is to consider partially hyperbolic (PH) ones~\cite{Hasselblatt2006}. The skew products over hyperbolic sets provide important examples of PH dynamics. Their central bundles are
tangent to the fibers.
What makes it worthy to consider the PH skew products is that by Hirsch--Pugh--Shub~\cite{HPS1977} theory, the systems \emph{conjugated} to such skew products form an open subset in PH. So such skew products are in a sense locally generic.

The Markov encoding of the base reduces them to the skew products over subshifts of finite type, see Definition~\ref{d:skpr}. The simplest skew products which correspond to $\dim E^c = 1$ have one-dimensional fibers. There is the single compact one-dimensional manifold without boundary, the circle~$S^1$, and the single one with boundary, the interval~$I = [0,1]$. As we already mentioned, the present work deals with skew products whose fibers are unit intervals.

%

Skew products are also widely used as a tool for the construction of (robust) examples of complicated behavior: convergence of orbits coexisting with minimality~\cite{Antonov1984, Kleptsyn2004}, non-removable zero Lyapunov exponents~\cite{GIKN2005}, attractors with intermingled basins~\cite{Kan1994, Bonifant2008, Ilyashenko2008}, multidimensional robust non-hyperbolic attractors~\cite{Viana1997}, and others. So it is very desirable to obtain a description of their dynamics.

To date, more attention has been paid to skew products with circle fibers than to those with interval fibers. The skew products with circle fibers turned out to exhibit some effects that cannot be observed for a single generic circle diffeomorphism. In his paper~\cite{Antonov1984}, Antonov proved that in an open set of step skew products (see Definition~\ref{d:step_and_mild}) with circle fibers the fiberwise coordinates of almost all orbits with the same base coordinate approach each other. Later in~\cite{Kleptsyn2004}, this effect was re-discovered by Kleptsyn and Nalski in the following terms. A generic step skew product with circle fibers may have a measurable ``attracting'' section, which is fiberwise approached by orbits of almost every point. At the same time the section is dense in the phase space. This construction was generalized to the case of mild (not necessarily step) skew products by Homburg in~\cite{Homburg2010}. Recall that the orbits of a single circle diffeomorphism either are dense (irrational rotation number) or tend to attracting periodic points (rational rotation number). The examples~\cite{Antonov1984},~\cite{Kleptsyn2004},~\cite{Homburg2010} exhibit a mixture of these two types of behavior.
%

The dynamical properties of generic diffeomorphisms of closed manifolds are known to be different from those of generic dissipative diffeomorphisms (that is, diffeomorphisms onto the image of a compact manifold with boundary) in the same dimension. For instance, a generic diffeomorphism of an interval has a finite number of fixed points. Any other point tends toward a fixed point. On the other hand, minimal diffeomorphisms of a circle are (metrically) generic.
This difference in behavior motivates us to go from studying skew products with circle fibers
to skew products with interval fibers.

Note that each skew product with interval fibers can be extended to a skew product with circle fibers. Namely, think of each interval fiber as of an arc of a circle and extend fiber maps to be circle diffeomorphisms. The simplest way to do it is to put a repelling fixed point at the complement to the arc. The resulting skew product with circle fibers can be made the same regularity of the fiber maps and the dependency of the fibers on the base as the initial interval skew product. Thus our result also describes the dynamics of an open (not dense) set of circle skew products.

In fact, our results are valid not only for the Markov measures but for a wider class of shift-invariant measures. The only condition we really need is stated in Proposition~\ref{p:measure_property}. Roughly speaking, it says that for any fixed ``right tail'' $(\om_1, \dots, \om_n, \dots)$ of the symbolic sequence in the base, the conditional probability to see any allowed symbol $p_j$ at the position~$\om_0$ is bounded from zero.

Similar statements also hold for generic \emph{mild} (see Definition~\ref{d:step_and_mild}) skew products with interval fibers; however, the genericity conditions must be changed. The proofs for this case require ideas beyond the present work, and we intend to present them in a separate paper.


\subsection{Remark}

In this paper we always assume that the fiberwise maps preserve the orientation of the unit interval. When this is not the case, our results can be derived as follows. Pass to the orientation-preserving covering. Namely, take two copies of the set of Markov states, marked with \textbf{+} and \textbf{-}. The orientation-preserving fiber maps take \textbf{+} to \textbf{+} and \textbf{-} to \textbf{-}. The orientation-changing fiber maps take \textbf{+} to \textbf{-} and \textbf{-} to \textbf{+}. Then our main result can be applied. After that we project the attractors and repellers back to the initial space.


\subsection{Outline of the paper}
The proof of our main result, Theorem~\ref{t:random}, is based on the link between the dynamics of skew products over one- and two-sided Markov shifts and the corresponding random dynamics on the interval.

In Sections~\ref{s:markov_def} and~\ref{s:one-sided} we recall some standard definitions, introduce the geometrical structures we need later, and state our main result for skew products, Theorem~\ref{t:markov_classification}. In Section~\ref{s:random_markov} we define random dynamics associated with step skew products. At the end of this section, we state the main result of the paper, Theorem~\ref{t:random}. This theorem describes the stationary measures of such random dynamics. In these two sections, an experienced reader can skip everything but the statements of the theorems.

In Section~\ref{s:typicalness} we give genericity conditions for Theorem~\ref{t:markov_classification} and Theorem~\ref{t:random}. Sections~\ref{s:proof} and~\ref{s:bax} are devoted to the proof of Theorem~\ref{t:random}: as we have already mentioned, this proof is based on the link between the skew products and the corresponding random dynamics (cf. Lemmas~\ref{l:absorb},~\ref{l:MFixed},~\ref{l:Markov-intervals},~\ref{l:Markov-uniqueness}). Finally, in Section~\ref{s:proof_markov_classif} we deduce Theorem~\ref{t:markov_classification} from Theorem~\ref{t:random}.

\section{Definitions and \hyperref[t:markov_classification]{the main theorem about  skew products}}\label{s:markov_def}

Suppose~$\sigma \colon \Sigma \to \Sigma$ is a transitive subshift of finite type (a topological Markov chain) with a finite set of states~$\{1, \ldots, N\}$ and $A = (a_{ij})_{i,j=1}^{N}$ is the transition matrix of $\sigma$, where $a_{ij} \in \{0,1\}$. Recall that $\Sigma$ is the set of all bilateral sequences~$\omega = (\omega_n)_{-\infty}^{+\infty}$ composed of symbols~$1,\ldots,N$ such that~$a_{\om_n \om_{n+1}} = 1$ for any~$n\in\bbZ$ (see, for instance,~\cite{Katok1995}).

The map~$\sigma$ shifts any sequence~$\om$ one step to the left: $(\sigma\omega)_n = \omega_{n+1}$.
By definition (pf.~\cite{Katok1995}), subshift is transitive iff
$$
\exists n \in \bbN \; \forall i,j \; (A^n)_{ij} > 0.
$$
Transitivity implies the indecomposability of the subshift. Indeed, for any~$m > 0$ the subshift~$\sigma^m$ with the same states allows one to go from any state to any other in finitely many steps. Thus for any~$m > 0$ the subshift~$\sigma^m$ cannot be split into two nontrivial subshifts of finite type.

We endow~$\Sigma$ with a metric defined by the formula
\begin{equation}    \label{e:markov_metric}
d (\om^1, \om^2) = \begin{cases}
2^{-\min\{|n| \, : \, \om_n^1\neq \om_n^2 \}}, &  \om^1 \neq \om^2, \\
0, &   \om^1 = \om^2,
\end{cases}
\qquad \om^1,\om^2 \in \Sigma.
\end{equation}

Now let~$M$ be a smooth manifold with boundary. Denote by~$\DiffIm^r(M)$ the space of $C^r$-smooth maps from $M$ to itself which are diffeomorphisms to their images.

\begin{Def} \label{d:skpr}
A \emph{skew product} over a subshift of finite type~$(\Sigma, \sigma)$ is a dynamical system~$F \colon \Sigma \times M \to \Sigma \times M$ of the form
$$
(\om, x) \mapsto (\sigma \om, f_{\om} (x)),
$$
where~$\om\in\Sigma$, $x\in M$, and the map~$f_{\om}(x) \in \DiffIm^1(M)$ is continuous in~$\om$. The phase space of the subshift is called the \emph{base} of the skew product, the manifold~$M$ is called the \emph{fiber}, and the maps~$f_{\om}$ are called the \emph{fiber maps}. \emph{The fiber over~$\omega$} is the set~$M_\om := \{\omega\} \times M \subset \Sigma\times M$.
\end{Def}

In any argument about the geometry of the skew products we always assume that the base factor of $\Sigma \times M$ is ``horizontal'' and the fiber factor is ``vertical''.

\begin{Def} \label{d:step_and_mild}
A skew product over a subshift of finite type is a \emph{step skew product} if the fiber maps~$f_{\om}$ depend only on the position~$\om_0$ in the sequence~$\om$. For the general skew products (fiber maps depend on the whole sequence~$\om$) we sometimes use the word \emph{mild}.
\end{Def}

The fiber~$M$ is always the unit interval~$I = [0,1]$ and the maps $f_\om \colon I \to f_\om(I)$
take the interval strictly inside itself. They also preserve the orientation, i.e., are strictly increasing. Suppose~$\sP$ is the set of all such skew products and~$\mcS \subset \sP$ is the subset of all step skew products. Note that~$\mcS$ is the Cartesian product of $N$ copies of $\DiffIm^1(I)$. We endow~$\sP$ with the metric
$$
\dist_\sP (F,G) := \sup_\om (\dist_{C^1} (f_\om^{\pm 1}, g_\om^{\pm 1})).
$$
This induces the $\max$-metric of a product on~$\mcS$.
In this paper, we consider only step skew products. We describe the dynamics of a generic~$F \in \mcS$.

In the rest of this text, all measures are assumed to be probabilities.
Like any dynamical system on a compact metric space, the subshift~$\sigma$ has a non-empty set of invariant measures. There is a natural class of them called \emph{Markov measures}. They are defined as follows.

Let~$\Pi = (\pi_{ij})_{i,j = 1}^N$, $\pi_{ij} \in [0,1]$ be a right stochastic matrix (i.e., $\forall i$ $\sum_j \pi_{ij} = 1$) such that $\pi_{ij} = 0$ iff~$a_{ij} = 0$. Let~$p$ be its eigenvector with non-negative components that corresponds to the eigenvalue~$1$:
\begin{equation}    \label{e:vector_p}
\forall i \; p_i \ge 0, \text{ and } \sum_i \pi_{ij} p_i = p_j.
\end{equation}
We can always assume~$\sum_i p_i = 1$.

For any finite word~$\om_k \ldots \om_m$, $k,m \in \bbZ$, $k \le m$ we consider a \emph{cylinder}
$$
C_w := \{ \om' \in \Sigma \mid \om_k' = \om_k, \ldots, \om_m' = \om_m \}.
$$
The cylinders form a countable base of the topology on~$\Sigma$. Thus a Borel measure on~$\Sigma$ is properly defined by its values on every cylinder.

\begin{Def}
$\nu$ is a \emph{Markov measure} constructed from the distribution~$p_i$ and the transition probabilities~$\pi_{ij}$ if its values on the cylinders are:
\begin{equation}    \label{e:markov_measure}
\nu(C_w) := p_{\omega_{k}} \cdot \prod\limits_{i=k}^{m-1} \pi_{\omega_{i} \omega_{i+1}}.
\end{equation}
\end{Def}

It is easy to see that the formula~\eqref{e:markov_measure} is consistent on the set of all cylinders. Thus~$\nu$ is well-defined. Moreover, it is invariant under the shift map~$\sigma$. Note that for any stochastic matrix~$\Pi$ there exists at least one vector~$p$ satisfying~\eqref{e:vector_p}. Such a vector is unique whenever the subshift is transitive and $\pi_{ij} \ne 0 \Leftrightarrow a_{ij} \ne 0$ (as in our case). If that is so, the measure~$\nu$ is ergodic; $\supp \nu$ coincides with~$\Sigma$.

Let~$\nu$ be any ergodic Markov measure on~$\Sigma$. From now on, the measure~$\nu$ is fixed.
\begin{Def} \label{d:std}
The \emph{standard measure}~$\mathbf{s}$ on $\Sigma \times I$ is the product of~$\nu$ and the Lebesgue measure on the fiber.
\end{Def}

Our goal is to give the description of the behavior of almost every orbit w.r.t. the standard measure. Such a description is given by Theorem~\ref{t:markov_classification} below. This is the main result of this paper.

\begin{Def}
We say that a closed set in the skew product is a \emph{bony graph} if it
intersects almost every fiber (w.r.t.~$\nu$) at a single point, and any
other fiber at an interval (a ``bone'').
\end{Def}

The name comes from the following simple observation. Any bony graph can be represented as a disjoint union of two sets, $K$ and $\Gamma$, where $K$ denotes the union of the bones. The projection of $K$ by $h$ to $\Sigma$ has zero measure, while $\Gamma$ is the graph of some measurable function $\varphi:\Sigma\setminus h(K)\to I$. By Fubini's Theorem, the standard measure of a bony graph is zero.

For any set~$B$, denote~$B_\om := B \cap I_\om$.
\begin{Def} \label{d:continuous-bony}
A bony graph~$B$ is a \emph{continuous-bony graph (CBG)} if $B_\om$ is upper-semicontinuous:
$$
\forall \om \, \forall \eps > 0 \quad \exists \delta > 0 \quad \text{ such that } \quad \dist(\om, \om') < \delta \Rightarrow B_{\om'} \subset U_\eps (B_\om).
$$
\end{Def}

In particular, the graph part~$\Gamma$ is a graph of a function which is continuous on its domain. 

In Theorem~\ref{t:markov_classification}, we show that the attractors and the repellers of generic skew products are CBGs. They are also $F$-invariant. It is easy to see that the restriction of the dynamics to an invariant CBG is very much the same as the restriction to the base:

\begin{Prop} \label{p:thin_conj}
Let $B$ be an invariant CBG.
Then the vertical projection~$h \colon B\setminus K \to \Sigma \setminus h(K)$ is a homeomorphism conjugating the dynamics of~$F|_{B\setminus K}$ and~$\sigma|_{\Sigma\setminus h(K)}$.
\end{Prop}

The following geometrical notions are crucial for our reasoning.
Let $\varphi_i \colon \Sigma \to I$ be two arbitrary functions, $\Gamma(\varphi_i)$ be their graphs, $i = 1,2$.

\begin{Def} \label{def:compare}
We write $\varphi_1 < \varphi_2$ if for any $\om \in \Sigma$
$$
\quad \varphi_1(\om) < \varphi_2 (\om).
$$
We also write~$\Gamma_{\varphi_1} < \Gamma_{\varphi_2}$ in this case.
\end{Def}

This definition admits a natural extension to the case of bony graphs:
\begin{Def} \label{def:compare-bony}
Let~$B_1, B_2$ be two bony graphs. We write $B_1 < B_2$ if for any $(\om, x_1) \in B_1$ and $(\om, x_2) \in B_2$ we have $x_1 < x_2$.
\end{Def}

Recall that a skew product permutes the fibers. Thus the image~$F(\Gamma)$ of any graph~$\Gamma$ is also a graph of some function.

\begin{Def} \label{def:drift_up}
We say that a graph~$\Gamma$ \emph{drifts up (down)} if $F(\Gamma) > \Gamma$ (respectively, $F(\Gamma) < \Gamma$).
\end{Def}

\begin{Def} \label{def:strip}
Let~$\varphi_1, \varphi_2 \colon \Sigma \to I$ be continuous functions such that~$\varphi_1 < \varphi_2$. The set
$$
S_{\varphi_1, \varphi_2} := \{ (\om,x) \mid \varphi_1(\om) \le x \le \varphi_2(\om) \}
$$
is called the \emph{strip between the graphs of~$\varphi_1$ and $\varphi_2$}.
\end{Def}

\begin{Def}
The strip~$S_{\varphi_1, \varphi_2}$ is \emph{trapping (nonstrictly trapping)} if $F(S_{\varphi_1, \varphi_2}) \subset \inter S_{\varphi_1, \varphi_2}$ (respectively $F(S_{\varphi_1, \varphi_2}) \subset S_{\varphi_1, \varphi_2})$. The strip is \emph{inverse trapping} (\emph{nonstrictly inverse trapping}) if the same holds true for $F^{-1}$.
%
\end{Def}

\begin{Rem} \label{r:monot_move}
Because the fiber maps are monotonous, $\varphi_1 < \varphi_2$ implies the inequality~$F^n(\varphi_1) < F^n (\varphi_2)$ for any~$n \ge 0$. This is also true for any~$n < 0$, provided that the preimages of the graphs are well-defined. Thus for any~$n \ge 0$ the image~$F^n(S_{\varphi_1, \varphi_2})$ is also a (non-empty) strip.
In particular, because $F(S_{\varphi_1, \varphi_2}) \subset S_{\varphi_1, \varphi_2}$, any trapping strip has a non-empty \emph{maximal attractor}:
\begin{equation}    \label{e:Amax}
A_{\textup{max}}(S_{\varphi_1, \varphi_2}) := \bigcap_{n=0}^{+\infty} F^n (S_{\varphi_1, \varphi_2}).
\end{equation}
\end{Rem}
In Theorem~\ref{t:markov_classification} we show that for a generic skew product the maximal attractor of any trapping strip is a CBG, provided the strip is indecomposable in a certain sense .

Among all invariant measures of a smooth dynamical system, the \emph{physical}, or \emph{SRB} (for Sinai, Ruelle, and Bowen) measures are of particular interest. See, for instance, the handbook~\cite{Barreira2006}. The following definition dates back to the classical papers~\cite{Sinai1972, Ruelle1976, Bowen2008}:
\begin{Def}
Let~$\mathbf{m}$ be an $F$-invariant measure. Consider the set~$V$ of all points~$p \in X$ such that for any function~$\varphi \in C(\Sigma\times I)$ the time average is equal to the space average:
$$
\lim\limits_{n\to\infty} \frac{1}{n}\sum_{i=0}^{n-1} \varphi (F^i(p)) = \int\limits_{\Sigma\times I} \varphi \, d\mathbf{m}.
$$
$\mathbf{m}$ is a \emph{physical measure} if the set~$V$ has a positive standard measure. The set~$V$ is then called its \emph{basin}.
\end{Def}

It turns out that in our case the following stronger property holds. The basin of any physical measure is a full-measure subset of an open subset of~$X$.

In general, very few invariant measures are physical. A long-standing open question asks if the attractors of typical dynamical systems carry physical measures (see, for instance, the review~\cite{Palis2005}). In his paper~\cite{Ruelle1976}, Ruelle proved that uniformly hyperbolic attractors do. In our case the attractors are partially hyperbolic instead, but we also manage to prove the existence of physical measures for them.

The main result of the paper is the following


\begin{Thm}\label{t:markov_classification}
For a generic~$F \in \mcS$ there is a finite collection of trapping strips and inverse trapping strips such that
\begin{enumerate}
\item\label{t:markov_classification:0.5} their union is the whole phase space~$X$;
\item\label{t:markov_classification:1} the maximal attractors~\eqref{e:Amax} of the trapping strips are CBGs; the repellers (i.e., the maximal attractors for the inverse map) of every inverse trapping strip are also CBGs;
\item\label{t:markov_classification:2} any trapping strip and any inverse trapping strip has a unique ergodic invariant measure such that its projection to the base is the Markov measure~$\nu$. Namely, it is the lift of~$\nu$ onto the attractor (or repeller) considered as the graph of a function defined almost everywhere. This measure is physical. Its basin contains a subset of full measure of the strip;
\item\label{t:markov_classification:2.5} for any trapping strip, there exists an invariant subset~$\Gamma$ of full measure (w.r.t. the measure from~\ref{t:markov_classification:2}) of its maximal attractor such that~$h\Gamma \subset \Sigma$ is a invariant residual subset of full measure~$\nu$. The vertical projection~$h$ conjugates the restrictions $F|_\Gamma$ and $\sigma|_{h\Gamma}$.
    Note that $h\Gamma$ is not closed but is invariant;
\item\label{t:markov_classification:3} the fiber-wise Lyapunov exponents of the attractor and repeller measures from~\ref{t:markov_classification:2} are non-zero;
\item\label{t:markov_classification:4} the graphs of the attractors and repellers are mutually comparable in the sense of Definition~\ref{def:compare}. Moreover, the attractors and the repellers are alternating: $\dots < A_1 < R_1 < A_2 < R_2 < \dots$.
\end{enumerate}
\end{Thm}

Here the set of all generic skew products is an open and dense subset of the space of all step skew products.
In~Section~\ref{s:typicalness} we explicitly state the genericity conditions.

We also provide an upper bound for the number of attractors and repellers, see Proposition~\ref{p:att_num}.

\begin{Rem} \label{r:finite_step}
One may consider the following generalization of step skew products. Let the fiber maps~$f_{\om}$ depend on finitely many symbols~$(\om_{-k} \ldots \om_{l})$, $k,l \in \bbN$, rather than on a single symbol. We call such systems \emph{multistep skew products}.

The statements of Theorem~\ref{t:markov_classification} also hold for generic multistep skew products. Indeed,
any multistep skew product $F$ can be replaced by the step skew product $G$ which base Markov states are the admissible words of length~$k+l+1$ of the base Markov chain of~$F$.
Note that the base Markov chain of $G$ is transitive if $F$'s one is.
\end{Rem}



\section{Skew products over one-sided shifts}   \label{s:one-sided}

%
%

Let~$\Sigma_+$ be the space of unilateral (infinite to the right) sequences $\omega=(\omega_n)_{0}^{+\infty}$ satisfying $a_{\om_n\om_{n+1}}=1$ for all~$n$. The left shift $\sigma_+ \colon \Sigma_+\to\Sigma_+, \;(\sigma_+ \om)_n=\om_{n+1}$ defines a non-invertible dynamical system on~$\Sigma_+$. The system $(\Sigma_+,\sigma_+)$ is a factor of the system~$(\Sigma,\sigma)$ under the ``forgetting the past'' map $\pi \colon (\om_n)_{-\infty}^{+\infty}\mapsto (\om_n)_{0}^{+\infty}$:
\begin{equation}\label{e:factor}
\begin{CD}
\Sigma @>{\sigma}>> \Sigma \\
@V{\pi}VV            @VV{\pi}V \\
\Sigma_+ @>{\sigma_+}>> \Sigma_+
\end{CD}
\end{equation}

Formulae~\eqref{e:markov_metric} and~\eqref{e:markov_measure} define the metric and the invariant measure~$\nu_+$ on~$\Sigma_+$. A measure is called \emph{invariant} under a non-invertible map~$F_+$ if for any measurable set~$A$
$$
\nu_+ (F_+^{-1} (A)) = \nu_+ (A).
$$

Recall that in a step skew product over the two-sided Markov shift the fiberwise maps depend only on $\om_0$. Thus one can pass from the skew product to the quotient:
\begin{equation}\label{e:skew-factor}
\begin{CD}
\Sigma \times I @>{F}>> \Sigma \times I \\
@V{\pi\times\Id}VV            @VV{\pi\times\Id}V \\
\Sigma_+ \times I @>{F_+}>> \Sigma_+ \times I,
\end{CD}
\end{equation}
where $F_+(\om,x)=(\sigma_+\om, f_{\om_0}(x))$.

There is a relation between the invariant measures of $F$ and $F_+$.
\begin{Prop}    \label{prop:meas_1_to_2}
\begin{enumerate}
\item For any~$F$-invariant measure~$\mathbf{m}$ its projection~$\mathbf{m}_+ = \pi_* \mathbf{m}$,
$$
\mathbf{m}_+(A) := \mathbf{m}(\pi^{-1} (A)),
$$
is $F_+$-invariant;
\item For any~$F_+$-invariant measure~$\mathbf{m}_+$ there exists an~$F$-invariant measure~$\mathbf{m}$ such that $\mathbf{m}_+ = \pi_* \mathbf{m}$.
\end{enumerate}
\end{Prop}

\begin{proof}
The first statement follows immediately from~\eqref{e:skew-factor}. To prove the second one, take an arbitrary measure~$\bar{\mathbf{m}}$ on $\Sigma \times I$ such that $\mathbf{m}_+ = \pi_* \bar{\mathbf{m}}$. Applying the Krylov-Bogolyubov averaging procedure~\cite{Krylov1947} to $\bar{\mathbf{m}}$, we obtain an~$F$-invariant measure~$\mathbf{m}$. And as the measure~$\mathbf{m}_+$ is $F_+$-invariant,~\eqref{e:skew-factor} implies that the projection of $\mathbf{m}$ is also equal to~$\mathbf{m}_+$.
\end{proof}


\section{Stationary measures of random walks}   \label{s:random_markov}

We assign a random walk to any one-sided skew product, as follows. Imagine that we are tracking the $F_+$-iterations of a point $(\om,x)\in\Sigma_+\times I$, where $\om$ is chosen randomly with respect to the measure~$\nu$. But we can only observe its $I$-coordinate and the symbol~$\om_0$. Then the sequence of our observations is a realization of a discrete Markov process~$\Pi(F)$ on the space $\mcI = \{1,\dots,N\}\times I$. The transition probability from a point~$(i,x)$ to a point~$(j,f_i(x))$ equals~$\pi_{ij}$.

For any measure~$\mu$ on the space~$\mcI$, it is natural to denote its stochastic image~$f_* \mu$ as
\begin{equation}\label{e:stoch_img}
  (f_* \mu)_j := \sum_i \pi_{ij} \cdot (f_j)_* \mu_i,
\end{equation}
where~$\mu_i$ is the restriction of the measure~$\mu$ to the interval~$I_i = \{i\} \times I$.

\begin{Def} \label{d:stationary}
A measure~$\mu$ on the space~$\mcI$ is \emph{stationary} if $f_* \mu = \mu$.
\end{Def}

Similarly to the invariant measures of dynamical systems, any stationary measure of a random processes admits a decomposition into an integral over \emph{ergodic} stationary measures, thus they are most important to study.

Any stationary measure of the process~$\Pi(F)$ induces an invariant measure of~$F_+$ that projects to the measure~$\nu_+$ on the base. Namely, a stationary measure~$\mu$ corresponds to the measure
\begin{equation}\label{e:inv_measure}
\mathbf{m}_+ (\mu) := \sum_k \nu_k^+ \times \mu_k,
\end{equation}
where $\nu_k^+$ is the conditional Markov measure on the cylinder $C_{+,k} = \{\, \omega \mid \omega_0 = k \} \subset \Sigma_+$ defined by~$\nu_+$:
$$
\nu_k^+=\frac{\nu_+|_{C_{+,k}}}{\nu_+(C_{+,k})}.
$$

On the other hand, Proposition~\ref{prop:meas_1_to_2} establishes a relation between invariant measures of $F_+$ and $F$. 
In Proposition~\ref{p:meas_on_limit} below we show that the ``attractor'' measures (see claim~\ref{t:markov_classification:2} of Theorem~\ref{t:markov_classification}) are projected under $\pi$-factorization to the measures of type~\eqref{e:inv_measure} which correspond to the stationary measures of~$\Pi(F)$.
However, the ``repeller'' invariant measures are mapped under $\pi$-factorization to measures which are supported on CBGs themselves. So, to study these measures, we will pass from~$F$ to the inverse skew product~$F^{-1}$.

\begin{Def}
A pair $(k,m)$ is \emph{admissible} if $\pi_{km} \neq 0$. The corresponding maps~$f_{k}\colon I_k \to I_m$ are also called \emph{admissible}.
\end{Def}

\begin{Def}
A subset of~$\mcI$ is a \emph{domain} if it intersects each $I_k$ at a nonempty interval.
\end{Def}

\begin{Def}
A domain $D=\bigsqcup_k D_k \subset \mcI$ is \emph{trapping (nonstrictly trapping)} if
any admissible map takes it to its interior (respectively, to itself). In other words,
$$
\forall k,m: \, \pi_{km} > 0, \, f_k(D_k) \subset \inter D_m \quad ( \, f_k(D_k) \subset D_m \, ).
$$
\end{Def}

It easy to see that the following proposition holds.
\begin{Prop}    \label{p:absorb}
The following conditions are equivalent:
\begin{enumerate}
\item the domain $D=\bigsqcup_k D_k \subset \mcI$ is trapping (nonstrictly trapping);
\item the strip
$$
\tilde{D}_+=\bigsqcup_k C_{+,k} \times D_k \subset \Sigma_+ \times I
$$
is a trapping (nonstrictly trapping) region for the skew product~$F_+$;
\item the strip
\begin{equation}\label{eq:Ds}
\tilde{D}=\bigsqcup_k C_k \times D_k \subset \Sigma \times I
\end{equation}
is trapping (nonstrictly trapping) for~$F$.
\end{enumerate}
Here, as before, $C_{+,k} = \{\, \omega \mid \omega_0 = k \} \subset \Sigma_+$, and $C_k = \{\, \omega \mid \omega_0 = k \} \subset \Sigma$.
\end{Prop}

\begin{Def}
A finite collection of domains is \emph{vertically ordered} if they can be enumerated in such a way that for any interval $I_k$ and for any $i<j$ the intersection of~$I_k$ with the $i$th domain is situated below the intersection of~$I_k$ with the $j$th domain.
\end{Def}

\begin{Def} \label{def:perehod}
An \emph{admissible composition} is a map of the form
$$
f_{w_1 \dots w_n} := f_{w_n}\circ \dots \circ f_{w_1}\colon I_{w_1}\to I_{w_{n+1}},
$$
where any pair of consequent symbols $(w_i,w_{i+1})$ is admissible. An admissible composition is called a \emph{simple transition} if all the symbols $w_i$, $i=1,\dots,n+1$ are different. An admissible composition is called a \emph{simple return} if $f_{w_1 \dots w_{n-1}}$ is a simple transition and $w_1=w_{n+1}$.
\end{Def}

For an arbitrary ergodic stationary measure~$\mu$ of the random process~$\Pi(F)$, denote by $\mu_k$ the restriction $\mu|_{I_k}$. Let $I_{\mu,k}=[A_{\mu,k}, B_{\mu,k}]$ be the interval that spans the support of $\mu_k$:
$$
A_{\mu,k}:= \min\supp \mu_k, \quad B_{\mu,k}:= \max\supp \mu_k.
$$
Note that the interval~$I_{\mu,k}$ may not coincide with~$\supp \mu_k$.

Finally, we state our main result on the stationary measures of~$\Pi(F)$:
\begin{Thm}\label{t:random}
Suppose~$F$ is a generic (in the sense of Section~\ref{s:typicalness}) step skew product. Then the following statements hold for the random process~$\Pi(F)$:
\begin{enumerate}
\item\label{t:random:1} there exist only finitely many ergodic stationary measures;
\item\label{t:random:2} the supports of these measures are contained in disjoint vertically ordered trapping domains;
\item\label{i:MFixed} for any ergodic stationary measure~$\mu$ and any $k$ there exist an attracting fixed point~$A$ of a simple return and a simple transition~$f$ such that~$A_{\mu,k} = f (A)$; the same is true for~$B_{\mu,k}$.
\item\label{t:random:3} the (random) Lyapunov exponents of these measures are negative.
\end{enumerate}
\end{Thm}

\begin{Rem}
An upper bound on the number of ergodic stationary measures is given in Corollary~\ref{c:att_num}.
\end{Rem}

\section{Genericity conditions}  \label{s:typicalness}

In this Section we explicitly state the genericity conditions for Theorem~\ref{t:markov_classification} and Theorem~\ref{t:random}:

\begin{enumerate}
\item\label{step_typicalnes:1} Any fixed point~$p$ of any simple return~$g$ is hyperbolic: $g'(p) \neq 1$;
\item\label{step_typicalnes:2} No attracting fixed point of a simple return is mapped to a repelling fixed point of a simple return by a simple transition. Also, no repelling fixed point of a simple return is mapped to an attracting fixed point of a simple return by a simple transition.
\item\label{step_typicalnes:3} One cannot choose from each interval $I_k$ a single point $a_k$ such that for any admissible couple~$(i,j)$ one would have $f_i(a_i)=a_j$ (in other words, so that the set~$\{a_k\}_{k=1}^{N}$ would be invariant under the random dynamics~$\Pi(F)$).
\end{enumerate}

\begin{Rem}
As there is only a finite number of simple transitions and simple returns, the above genericity conditions are satisfied on the complement to a finite union of codimension one subsets of~$\mcS$.
\end{Rem}
\begin{Rem}
The condition~\ref{step_typicalnes:3} is closely related to the \emph{accessibility} property of partially hyperbolic systems (see, for instance, the handbook~\cite{Barreira2006}).
\end{Rem}
\begin{Rem}
Instead of two latter conditions one may impose the following: there is no way to take a fixed point of a simple return to a fixed point of a simple return by a simple transition map, except for the case when the word defining the simple transition is a suffix of the word defining the latter simple return.
\end{Rem}

\begin{Rem}
  In general, Theorems~\ref{t:random} and~\ref{t:markov_classification} are false without these assumptions. For example, Ilyashenko in~\cite{Ilyashenko2010b} gave an example of ``thick'' (positive measure) attractors appearing when the condition~\ref{step_typicalnes:2} is broken.
\end{Rem}

\section{Proof of the \hyperref[t:random]{stationary measures theorem}}  \label{s:proof}

Let~$\mu$ be an arbitrary ergodic stationary measure for the random process~$\Pi(F)$ described in Section~\ref{s:random_markov}. For any admissible transition~$(i,j)$ we have
$$
f_{i} (\supp \mu_i) \subset \supp\mu_j.
$$
Because the maps~$f_k$ are monotonous, the disjoint union of the intervals~$I_{\mu,k}$ is forward-invariant: for any admissible transition~$(i,j)$
\begin{equation}    \label{e:I_k_forw_inv}
f_{i} (I_{\mu,i}) \subset I_{\mu,j}.
\end{equation}
Thus the domain~$\mcI_\mu = \bigsqcup_k I_{\mu,k}$ is nonstrictly trapping.

\begin{Rem}
The genericity conditions of Section~\ref{s:typicalness} imply that no interval~$I_{\mu,k}$ can be a single point. Otherwise by~\eqref{e:I_k_forw_inv} any interval~$I_{\mu,k}$ is a single a point. And this contradicts with the genericity condition~\ref{step_typicalnes:3}.
\end{Rem}


\begin{Lem}   \label{l:absorb}
There exist arbitrarily small strictly trapping neighborhoods of $\mcI_\mu$.
\end{Lem}

Before proving Lemma~\ref{l:absorb}, we need some auxiliary statements. The following Lemma is the statement~\ref{i:MFixed} of Theorem~\ref{t:random}.

\begin{Lem} \label{l:MFixed}
For any $k$ there exist an attracting fixed point~$A$ of a simple return and a simple transition~$f$ such that~$A_{\mu,k} = f (A)$.
The same is true for~$B_{\mu,k}$.
\end{Lem}

\begin{proof}
By the Kakutani random ergodic theorem~\cite{Kakutani1951} (see also~\cite{Furman2002}) a generic sequence of random iterations~$(k_n, x_n)$, $x_n \in I_{k_n}$ of a $\mu$-generic initial point is distributed with respect to the measure~$\mu$. Let us choose and fix such a generic initial point~$(k_0,x_0)$, different from $(k_0, A_{\mu,k_0})$ and $(k_0,B_{\mu,k_0})$. Denote by~$(k_n,x_n)$ a generic sequence of iterations of~$(k_0,x_0)$. The points~$(k_n,x_n)$ are distributed w.r.t.~$\mu$.

Thus for any $k$ the set~$X_k = \{x_n\}|_{k_n=k}$ is dense in~$\supp \mu_k$. In particular, the points~$(k,A_{\mu,k})$ and $(k,B_{\mu,k})$ are accumulation points of~$X_k$.

Now we modify the sequence~$(k_n,x_n)$ so that it becomes monotonous in the following sense:

\begin{Def}
An admissible sequence of iterations is called \emph{downwards monotonous} if for any $(k_m,x_m)$, $(k_n,x_n)$ such that $k_m = k_n$ and $m < n$ we have $x_m > x_n$. In the same way we define \emph{upwards monotonicity}.
\end{Def}

The following proposition selects a monotonous subsequence from any finite sequence of iterations.

\begin{Prop}\label{p:downwards}
For any finite admissible sequence of iterations
$$
(w_0,x_0),\dots,(w_n,x_n)
$$
there exists downwards monotonous finite admissible sequence
$$
(w_0,x_0)=(w'_0,x'_0),\dots,(w'_{n'},x'_{n'}),
$$
such that $w'_{n'} = w_n$ and $x'_{n'} \le x_n$.
\end{Prop}
\begin{proof}
To find such a sequence, let us write out the symbols of the original sequence, and remove each simple return such that the point in fiber after that return is higher than before.

More formally, the proof is by induction on~$n$. Assume that the existence of the desired word~$M(w)$ is proven for any initial word~$w$ of length~$|w| \le n$.
Let us now prove it for the word~$w_0 \ldots w_n$ of length $n+1$.

Denote $w = w_0 \dots w_{n-1}$. If the symbol~$w_n$ is not contained in~$M(w)$, then the word $w':=M(w)w_n$ is the desired one.
Otherwise consider the last occurrence of~$w_n$ in the word~$M(w)$: let
$M(w)=uw_n v$, where the symbol~$w_n$ does not appear in the word~$v$. Compare the images of the initial point under the action of the words~$uw_n$ and~$uw_n vw_n= M(w)w_n$ (both these images belong to the interval~$I_{w_n}$). If the second of these images lies below the first one, the conclusion of the lemma is satisfied for the word $w':=M(w) w_n$. Otherwise one can take~$w' := uw_n$.
\end{proof}

Recall that the original sequence of iterations of the point $(k_0,x_0)$ is dense in~$\supp\mu$. Combining this with the previous proposition, we obtain the following
\begin{Prop}    \label{p:word_constr}
There exists a sequence of words $w^j$, starting with the symbol~$k_0$ and ending with~$k$, such that for any of these words the corresponding iterations of the point $(k_0,x_0)$ are monotonous downwards, and the final images~$(k,x_j)$ tend to~$(k,A_{\mu,k})$ as $j \to \infty$. Also, the sequence of lengths~$|w^j|$ tends to infinity.
\end{Prop}

\begin{proof}
First, take a subsequence of iterations of~$(k_0,x_0)$ that tends to~$(k,A_{\mu,k})$. Then apply~Proposition~\ref{p:downwards} to obtain a sequence of words~$w_j$ such that the sequence~$(k, x_j)$ is monotone. Because $I_{\mu,k}$ is invariant, for every~$j$ we have $A_{\mu,k} \le x_j \le B_{\mu,k}$. By the sandwich rule, $(k,x_j) \to (k,A_{\mu,k})$ as $j \to +\infty$.

Now assume the lengths~$|w^j|$ do not tend to infinity.
Then there exists an admissible composition~$G$ that takes the initial point~$(k_0,x_0)$ exactly to the point~$(k,A_{\mu,k})$. But this is impossible. Indeed, otherwise the image of the point~$P = (k_0,A_{\mu,k_0})$ under $G$ must be below~$(k,A_{\mu,k})$, because the point $P$ is below~$(k_0,x_0)$. And at the same time $G(P) \in \supp \mu \subset [A_{\mu,k}, B_{\mu,k}]$.
\end{proof}

Now we can conclude the proof of Lemma~\ref{l:MFixed}. Indeed, consider the sequence of words~$(w^j)$ constructed in Proposition~\ref{p:word_constr}. As their lengths tend to infinity, at some point they exceed $N$; in particular, there are repeating symbols in these words.

Let $|w^j| > N$. Passing from the end of this word to the beginning, find the first repetition of symbols. Namely, let
$$
w^j=asbsc,
$$
where $a,b,c$ are words, $s$ is a symbol, and all the symbols in the word~$bsc$ are different.

So the composition $f_{sb}$ is a simple return. It takes the point~$f_{a}(x_0)$ to a point below~$f_{a}(x_0)$. Meanwhile, both these points belong to~$\supp\mu$. Thus for any $M \ge 0$ the point~$f_{sb}^M(f_{a}(x_0))$ also belongs to~$\supp\mu$. The sequence~$f_{sb}^M(f_{a}(x_0))$ tends to an attracting fixed point~$p^j$ of the map~$f_{sb}$ as $M \to +\infty$. Hence~$p^j\in \supp\mu$.

Denote by~$g^j$ the simple transition~$f_{sc}$. Note that~$f_{sc}(p^j) \le f_{asbsc}(x_0)$ because $p^j \le f_{asb}(x_0)$. So the sequence $f_{w^j}(x_0)$ majorizes the sequence $g^j(p^j)$. Hence the latter also converges to~$(k,A_{\mu,k})$.

But there is only a finite number of simple returns and of simple transitions. By the genericity conditions any simple return has finitely many fixed points. So there are only finitely many different maps~$g^j$ and points~$p^j$. The sequence we have just constructed ranges over a finite set. The fact this sequence converges to~$(k,A_{\mu,k})$ implies that the point~$(k,A_{\mu,k})$ belongs to this finite set, so $(k,A_{\mu,k}) = g^j (p^j)$ for some $j$.

The same reasoning proves the statement for~$B_{\mu,k}$.
\end{proof}

\begin{proof}[Proof of Lemma~\ref{l:absorb}]
First, make an observation that will be useful later in the proof.
Let~$G$ be any simple return. Then either~$A_{\mu,k}$ is a fixed point of~$G$ or $G(A_{\mu,k}) \in (A_{\mu,k}, B_{\mu,k})$. In the first case, $A_{\mu,k}$ is hyperbolic by genericity condition~\ref{step_typicalnes:1}, and it is attracting by Lemma~\ref{l:MFixed} and genericity condition~\ref{step_typicalnes:2}.
The same holds for the point~$B_{\mu,k}$.

Denote by~$\Phi(D)$ the \emph{diffusion} of a set $D \subset \mcI$:
$$
\Phi(D)= \bigsqcup_{m} \left(\bigcup_{k: \, (k,m) \, \text{admissible}} f_k(D\cap I_{k})\right),
$$
$k,m \in \{1, \dots, N\}$.
Note that a domain $D$ is non-strictly trapping if and only if $\Phi(D)\subset D$, and trapping iff $\Phi(D)\subset \inter D$. (Note also that the domain $D$ is nonstrictly expelling, that is, the strip corresponding to $D$ by~\eqref{eq:Ds}, is nonstrictly expelling, if and only if $\Phi(D)\supset D$. In particular, $\Phi(D)$ is not necessary a superset of~$D$.)

To construct the desired trapping domain, we will first prove the following statement. For any $\eps>0$,
consider closed $\eps$-neighborhood~$D=\bigsqcup_k \overline{U_{\varepsilon} (I_{\mu, k})}$ of the set~$\mcI_\mu$. Then, for any sufficiently small $\eps>0$ the $N+1$-th image of $D$ is inside the interior of the union of the first $N$ images:
\begin{equation}\label{eq:N+1}
\Phi^{N+1}(D) \subset \inter \bigcup_{j=0}^N \Phi^j(D).
\end{equation}

Indeed, take any admissible word $w=w_1\dots w_{N+1}$ of length~$N+1$. We want to study the image $f_{w_1\dots w_N}:I_{w_1}\to I_{w_{N+1}}$ of $[A_{\mu,w_1}-\eps,B_{\mu,w_1}+\eps]$; it suffices to show that both endpoints of this image belong to the interior of the right hand-side of~\eqref{eq:N+1}. Also, it suffices to consider only the images of $A_{\mu,w_1}-\eps$: for the second endpoint, the argument is analogous.

Consider first all the intermediate images of the non-shifted endpoint $f_{w_1\dots w_j}(A_{\mu,w_1})$. Note, that if at least one of them belongs to the interior of the corresponding $I_{\mu, w_{j+1}}$, then so does the final image $f_{w_1\dots w_N}(A_{\mu,w_1})$, and thus by continuity for any sufficiently small $\eps$ so does $f_{w_1\dots w_N}(A_{\mu,w_1}-\eps)$.

Otherwise, if all the images $f_{w_1\dots w_j}(A_{\mu,w_1})$ exactly coincide with the corresponding endpoints $A_{\mu,w_{j+1}}$, note that in the word $w$ there is at least one simple return $f_{w_i\dots w_{j-1}}$ (so that $w_i=w_j$). Let $G=f_{w_i\dots w_{j-1}}:I_{w_i}\to I_{w_i}$ be the corresponding return map. Then, the point $A_{\mu,w_i}$ is an attracting fixed point of $G$ (due to the argument in the beginning of the proof), and we have $A_{\mu,w_i}=f_{w_1\dots w_{i-1}}(A_{\mu,w_1})=f_{w_1\dots w_{j-1}}(A_{\mu,w_1})$. Hence, for any sufficiently small $\eps>0$,
$$
f_{w_1\dots w_{i-1}}(A_{\mu,w_1}-\eps)< G\circ f_{w_1\dots w_{i-1}}(A_{\mu,w_1}-\eps)  = f_{w_1\dots w_{j-1}}(A_{\mu,w_1}-\eps),
$$
what implies
$$
f_{w_1\dots w_{i-1}w_j\dots w_N}(A_{\mu,w_1}-\eps)<f_{w_1\dots w_{N}}(A_{\mu,w_1}-\eps).
$$
The point $f_{w_1\dots w_{N}}(A_{\mu,w_1}-\eps)$ belongs to the interior of $\cup_{j=0}^N \Phi^j(D)$.
So, for any admissible word $w$ of length $N+1$ we have the desired inclusion, and~\eqref{eq:N+1} is proven.


Now, denote by $\Phi_{\delta}$ the \emph{$\delta$-dispersed diffusion}:
$$
\Phi_\delta (D) := \overline{U_\delta (\Phi(D))}.
$$
For any sufficiently small $\eps>0$, as the inclusion~\eqref{eq:N+1} is stable under small perturbations, we have for any sufficiently small $\delta>0$
\begin{equation}\label{eq:delta}
\Phi_{\delta}^{N+1}(D) \subset 
\inter \bigcup_{j=0}^N \Phi_{\delta}^j(D).
\end{equation}

Consider then for any such $\eps$, $\delta$ the domain $\tilde{D}_{\delta}:=\bigcup_{j=0}^{N} \Phi_\delta^j (D)$. Immediately from definition we see that for any domain $Y$ one has $\Phi(Y)\subset \inter \Phi_{\delta}(Y)$, in particular,
$$
\Phi(\Phi_{\delta}^{j}(D))\subset \inter \Phi_{\delta}^{j+1}(D), \quad j=0,\dots N.
$$
Hence,
$$
\Phi(\tilde{D}_{\delta}) \subset \inter \tilde{D}_\delta,
$$
and $\tilde{D}_\delta$ is a trapping domain.
\end{proof}

Any ergodic stationary measure~$\mu$ is supported inside some non-strictly trapping domain that is a union of intervals~$\mcI_\mu = \bigsqcup_k I_{\mu,k}$. The ends of these intervals are images of fixed points of simple returns under simple transitions. So there are only finitely many possibilities for such trapping domains.
Now the conclusions~\ref{t:random:1}, \ref{t:random:2} of Theorem~\ref{t:random} are reduced to the following two lemmas.

\begin{Lem}\label{l:Markov-intervals}
For any two ergodic stationary measures~$\mu_1$ and~$\mu_2$ the corresponding intervals~$I_{\mu_1,k}$ and~$I_{\mu_2,k}$ are either disjoint for any $k$ or coincide for any $k$. In the former case they are situated in the same order on all the intervals~$I_{k}$.
\end{Lem}

\begin{Lem}\label{l:Markov-uniqueness}
Any stationary ergodic measure $\mu$ is the unique stationary measure in the corresponding trapping domain~$\mcI_\mu$. 
\end{Lem}

To prove them, we need the following (useful for many reasons) lemma and its corollaries.

\begin{Lem}\label{l:Markov-contractions}
Let~$\mu$ be an ergodic stationary measure, $k$ and $k'$ be two arbitrary symbols. Then for any~$\eps > 0$ there exist admissible compositions~$G_{A, B} \colon I_k \to I_{k'}$ such that
$$
G_A (I_{\mu,k}) \subset U_\eps (A_{\mu,k'}), \quad G_B (I_{\mu,k}) \subset U_\eps (B_{\mu,k'}).
$$
\end{Lem}

\begin{proof}
Apply the genericity condition~\ref{step_typicalnes:3} (sf. Section~\ref{s:typicalness}) to the set of lower endpoints~$\{A_{\mu,k}\}$ of intervals~$I_{\mu,k}$. There exist~$k''$ and an admissible composition~$G_1 \colon I_k \to I_{k''}$ such that $p := G_1(A_{\mu,k}) > A_{\mu,k''}$. Because $A_{\mu,k''} = \inf \supp \mu_{k''}$, $\mu([A_{\mu,k''},p]) > 0$. Thus the generic start point~$(x_0,k'')$ from the proof of Lemma~\ref{l:MFixed} can be chosen from~$[A_{\mu,k''},p]$.

Because~$\mu$ is ergodic, a random orbit with a $\mu$-generic start point is dense in $\supp\mu$ almost surely. So for any $\eps > 0$ there exists an admissible composition~$G_2 \colon I_{k''} \to I_{k'}$ such that~$G_2(x_0, k'') \in U_\eps (B_{\mu,k'})$.



The monotonicity of fiberwise maps~$f_m$ implies that the image~$(G_2 \circ G_1) (I_{\mu,k})$ is situated on the interval~$I_{k'}$ strictly above the point~$G_2(x_0, k'')$. The map~$G_B := G_2 \circ G_1$ is constructed. The procedure for the lower endpoint~$A_{\mu,k'}$ is analogous.
\end{proof}

\begin{Cor} \label{cor:nonintersect}
For any interval~$I_k$ there exist an interval~$I_m$ and admissible compositions $G_{1, 2} \colon I_m \to I_k$ such that the images~$G_1 (I_{\mu,m})$ and~$G_2(I_{\mu,m})$ do not intersect.
\end{Cor}

\begin{proof}
Indeed, it suffices to take the compositions that send the interval $I_k$ to disjoint neighborhoods of the points~$A_{\mu,m}$ and $B_{\mu,m}$ respectively.
\end{proof}

\begin{Cor} \label{cor:meas_not_coincide}
For any interval~$I_k$ there exist an interval~$I_m$ and an admissible composition $G \colon I_m \to I_k$ such that $\mu_m \neq G_* \mu_k$.
\end{Cor}

\begin{proof}[Proof of Lemma~\ref{l:Markov-intervals}]
Assume the converse: let the intervals $I_{\mu_1,k}$ and $I_{\mu_2,k}$ intersect but not coincide. Then there is an endpoint of one of them that does not belong to another; without loss of generality let it be the point~$A_{\mu_1,k}$. Take a neighborhood~$U \ni A_{\mu_1,k}$ such that $I_{\mu_2,k} \cap U = \emptyset$.

By Lemma~\ref{l:Markov-contractions} applied for $k=k'$, there exists an admissible return~$G$ such that $G(I_{\mu_1,k}) \subset U$. On the other hand, any interval~$I_{\mu,k}$ is absorbing under any admissible return, so~$G(I_{\mu_2,k}) \subset I_{\mu_2,k}$. Then the nonempty set~$G(I_{\mu_1,k} \cap I_{\mu_2,k})$ is contained in both~$U$ and~$I_{\mu_2,k}$. This contradiction proves the first part of the lemma.

The second part of Lemma~\ref{l:Markov-contractions} follows directly from the monotonicity of the fiberwise maps.
\end{proof}

To go on with the proof of Lemmas~\ref{l:Markov-intervals} and~\ref{l:Markov-uniqueness}, we temporary \textbf{switch to skew product mode:}

As the diffeomorphisms~$f_{\omega}$ take the interval~$I$ to its interior, the inverse maps~$f_{\omega}^{-1}$ are not everywhere defined. To overcome this technical obstacle, we extend the interval fibers~$I$ to circle fibers~$S^1 \supset I$ in the following way:
\begin{itemize}
\item each map~$f_\om$ is an orientation-preserving diffeomorphism of the circle;
\item on the set~$S^1 \setminus I$ each map~$f_\om$ is strictly stretching: $f_\om' > 1$;
\item the whole skew product is generic in the sense of Section~\ref{s:typicalness}.
\end{itemize}

Then for the inverse skew product there exists a single attractor~$J$ in the strip~$\Sigma \times (S^1 \setminus I)$. Moreover, it is a (non-bony) graph of a continuous function~$j \colon \Sigma \to S^1$. This statement is straightforward. The detailed proof is given, for instance, in~\cite[Props.~2, 3]{Ilyashenko2010}.

For any $n \in \bbZ$ denote by~$f_{n,\om} (x)$ the~$n$th iterate of a fiber point~$x$ with a sequence~$\om$ in the base, so
\begin{equation}    \label{e:n_image}
F^n(\om,x) = (\sigma^n \om, f_{n,\om} (x)).
\end{equation}
Then for any $n > 0$
$$
f_{n,\om} = f_{\sigma^{n-1}\om} \circ f_{\sigma^{n-2}\om} \circ \ldots \circ f_{\om},\;\text{ and }\;
f_{-n,\om} = f_{\sigma^{-n}\om}^{-1} \circ f_{\sigma^{-n+1}\om}^{-1} \circ \ldots \circ f_{\sigma^{-1}\om}^{-1} = f^{-1}_{n,\sigma^{-n}\om}.
$$
In particular, for step skew products
$$
f_{n,\om} = f_{\om_{n-1}} \circ f_{\om_{n-2}} \circ \ldots \circ f_{\om_0},\;\text{ and }\;
f_{-n,\om} = f_{\om_{-n}}^{-1} \circ f_{\om_{-n+1}}^{-1} \circ \ldots \circ f_{\om_{-1}}^{-1}.
$$

Suppose $\mu$ is an ergodic stationary measure of $\Pi(F)$. By Lemma~\ref{l:absorb}, the non-strictly trapping domain~$\mcI_\mu$ has arbitrary small strictly trapping neighborhood~$U(\mcI_\mu) = \sqcup_k U_k$. Denote by~$\tilde D \subset \Sigma \times S^1$ the corresponding strip in the two-sided skew product:
$$
\tilde D := \bigsqcup_k C_k \times U_k,
$$
where~$C_k = \{\, \omega \mid \omega_0 = k \} \subset \Sigma$. By Proposition~\ref{p:absorb}, $\tilde D$ is also strictly trapping.

\begin{Prop}\label{p:Markov-bony}
The maximal attractor~$\Amax$ of the trapping strip $\tilde D$ is an invariant continuous-bony graph.
\end{Prop}
\begin{proof}

Like before, we denote by~$A_\om$ the intersection of $\Amax$ and the fiber $\{ \om \} \times I$. By definition,
$$
A_\om = \bigcap_{n \ge 0} A_\om^{(n)},
$$
where
$$
A_\om^{(n)} = f_{-1} \circ \dots \circ f_{\om_{-n}} (\tilde D_{\sigma^{-n} \om}).
$$
Note that~$A_\om^{(n)}$ is a sequence of nested intervals, and thus $A_\om$ is either an interval or a single point. Also note that if some sequences $\om$ and $\om'$ are close enough to each other, say,
$$
\om_{-n}' = \om_{-n}, \dots, \om_{-1}' = \om_{-1},
$$
then $A_\om^{(n)} = A_{\om'}^{(n)} \supset A_{\om'}$.
This implies the upper-semicontinuity of $A_{\omega}$. This semi-continuity, once we prove that $A_{\max}$ is a bony graph, will immediately imply the continuity of its graph part. The semi-continuity also implies that the set $\{ \om\,|\, A_\om \text{ is a point}\}$ is residual.

To prove that~$\Amax$ is actually a bony graph, we use an argument similar to the~\cite[Theorem 3]{Kudryashov2010} by Kudryashov. By the Fubini Theorem, $\Amax$ is a bony graph whenever its standard measure~$\mathbf{s}$ (see Definition~\ref{d:std}) is zero. For every $k = 1 \dots N$ and $x \in I_{\mu,k}$, denote by~$\Om_{k,x} \subset \Sigma$ the slice of~$\Amax \cap \{ \om_0 = k\}$ by the horizontal line $\Sigma \times \{x\}$. Because~$F^{-1}(\Amax) = \Amax$, $\Om_{k,x}$ is the set of sequences such that
\begin{itemize}
  \item $\om_0 = k$;
  \item $\forall n \ge 0 \quad F^{-n} (\om,x) \in \tilde D$.
\end{itemize}
By definition, $\Amax = \sqcup_{k,x} \Om_{k,x} \times \{x\}$.

Now we show that for all $k$, for all $x \in I_{\mu,k}$ we have $\nu(\Om_{k,x}) = 0$.
By the Fubini Theorem, this will imply $\mathbf{s}(\Amax) = 0$. In a sense, we have just switched the order of integration for $\Amax$. 

To prove this, we cover~$\Om_{k,x}$ by a disjoint union of cylinders of arbitrary small $\nu$-measure.

\subsection*{1st generation of the cylinders}

By Corollary~\ref{cor:nonintersect} of Lemma~\ref{l:Markov-contractions}, there exist two words~$w^{(1)}, w^{(2)}$ of same length~$L$, $w^{(i)} = w^{(i)}_{-L} \dots w^{(i)}_{-1}$ such that for any $i=1, 2$
\begin{itemize}
  \item the word $w^{(i)} k$ is admissible;
  \item $w_{-L}^{(i)} = k$;
  \item $f_{w^{(1)}} (I_{\mu,k}) \cap f_{w^{(2)}} (I_{\mu,k}) = \emptyset$.
\end{itemize}
The latter implies that for at least one of these words, denote it by $w$, holds~$f_w^{-1} (x) \notin I_{\mu,k}$. Now for every $\om \in \Om_{k,x}$ we know that $\om_{[-L,-1]} \ne w$, so
$$
\Om_{k,x} \subset \bigsqcup_{s \ne w} \left\{ \om_{[-L,-1]} = s \right\}.
$$
The right-hand side is the 1st generation of the cylinders.

\subsection*{The next generation of the cylinders}

Now we subdivide each of the 1st generation cylinders and see that
$$
\left\{ \om_{[-L,-1]} = s \right\} = \bigsqcup_{s'} \left\{ \om_{[-2L,-1]} = s's \right\}.
$$
For each $s \ne w$, we now have new $x_s := f_s^{-1} (x) \in I_{\mu,k}$ and a new word $w' \in \{w^{(1)}, w^{(2)}\}$ such that $f_{w'}^{-1} (x_s) \notin I_{\mu,k}$.
Thus
\begin{equation}\label{e:next_cyl}
\Om_{k,x} \cap \left\{ \om_{[-L,-1]} = s \right\} \subset \bigsqcup_{s' \ne w'} \left\{ \om_{[-2L,-1]} = s's \right\}.
\end{equation}

\subsection*{Uniform decrease of the cylinders}

Let~$w$ be any admissible word. For any cylinder $C = \left\{ \om_{[n_1, n_2]} = u\right\}$, $u$ being some word, denote by~$wC$ the cylinder
$$
wC := \left\{ \om_{[n_1 - L, n_2]} = wu \right\},
$$
provided $wu$ is also admissible. In other words, we concatenate $w$ and $u$. Obviously, $wC \subset C$.

\begin{Prop}    \label{p:measure_property}
  There exists~$0 < \lambda < 1$ such that for any~$w \in \{ w^{(1)}, w^{(2)} \}$ and for any cylinder~$C$ such that~$wu$ is admissible, we have
  $$
  \frac{\nu(wC)}{\nu(C)} > \lambda.
  $$
\end{Prop}

This is an obvious corollary of $\nu$ being a Markov measure.
In fact, this proposition seems to hold for a much wider class of measures than the Markov ones. For instance, we suspect that it holds for all the measures which appear as SRB in smooth partially hyperbolic dynamics, and perhaps for all or a large class of Gibbs measures. It seems to be interesting on its own to determine the class of measures for which Proposition~\ref{p:measure_property} holds. Once one can prove Proposition~\ref{p:measure_property} for a larger class of measures, this will immediately generalize Theorem~\ref{t:markov_classification} to them.


\subsection*{End of the proof of Proposition~\ref{p:Markov-bony}}

On each step of the construction, we subdivide the cylinders in the same way. Because of Proposition~\ref{p:measure_property}, each next generation of cylinders has at least~$\lambda$ part of the ``bad'' ones, i.e. those who eventually throw $x$ out of $\tilde D$. So each subdivision reduces the $\nu$-measure of the cover by at least $1-\lambda < 1$ times. So we are able to cover~$\Om_{k,x}$ with a set of an arbitrary small $\nu$-measure. Thus~$\Om_{k,x} = 0$, and Proposition~\ref{p:measure_property} is proved.

\end{proof}

\begin{Prop}    \label{p:SRB}
The lift~$\mathbf{m}_{\varphi}$ of the measure~$\nu$ on the graph part~$\Gamma_{\varphi}$ of~$\Amax$ is physical.
\end{Prop}
\begin{proof}
Assume that the point~$(\om,x)\in S_{\varphi_1,\varphi_2}$ is such that the sequence of the time averages of its base coordinate~$\om\in\Sigma$ converges to the measure~$\nu$. (Because the measure~$\nu$ is ergodic, this holds for the set of points of full standard measure in~$S_{\varphi_1,\varphi_2}$.)

Consider the sequence of time averages
$$
\theta_n:=\frac{1}{n} \sum_{j=0}^{n-1} \delta_{F^j(\om,x)}.
$$
We have to show that this sequence converges (in $*$-weak topology) to the measure~$\mathbf{m}_{\varphi}$. To do so, note, on one hand, that any limit point $\theta=\lim_i \theta_{n_i}$ of this sequence of measures due to the choice of the coordinate~$\omega$ is projected to the measure $\nu$ on the base. Moreover, its support is contained in the maximal attractor $A_{\max}:=\cap_{j} F^{j}(S_{\varphi_1,\varphi_2})$. On the other hand, $\nu$-almost every fiber $h^{-1}(\omega')$ intersects the maximal attractor by a single point~$(\om', \varphi(\om'))$. Thus the conditional measure of~$\theta$ on almost every fiber is the Dirac measure~$\delta_{(\om', \varphi(\om'))}$, and hence the measure $\theta$ itself coincides with the measure~$\bf{m}_{\varphi}$.

The uniqueness of the limit point and compactness of the space of measures now imply that all the sequence of time averages~ $\theta_n$ converges to the measure~$\bf{m}_{\varphi}$.
\end{proof}
\begin{Rem}\label{r:lim}
The same arguments imply than for any measure $\textbf{m}$, supported in $S_{\varphi_1,\varphi_2}$, that projects to the Markov measure $\nu$ on the base, its iterations converge to the measure~$\mathbf{m}_{\varphi}$. Indeed, any accumulation point of the sequence of measures $F_*^n \mathbf{m}$ is a measure, supported on the maximal attractor $A_{max}$ and such that it projects to the Markov measure on the base. Hence, any accumulation point of this sequence of iterations coincides with $\mathbf{m}_{\varphi}$ and thus the entire sequence converges to~$\mathbf{m}_{\varphi}$.
\end{Rem}

Now to prove Lemma~\ref{l:Markov-uniqueness}, we \textbf{switch back
back to Random Walk mode:}
\begin{proof}[Proof of Lemma~\ref{l:Markov-uniqueness}]
%
As we have already mentioned in Section~\ref{s:random_markov}, in a skew product over the Markov shift in the space of one-sided sequences~$\Sigma_+$
to the measure~$\mu'$ corresponds measure~$\mathbf{m}_+'$, defined by formula~\eqref{e:inv_measure}.
As one can easily see from the definition, $\mathbf{m}_+'$ is $F_+$-invariant.

Consider now the measure~$\mathbf{m}' = \sum_k \nu_k \times \mu_k$, that is given by the same sum, but for the skew product over the shift on space of two-sided sequences. In general, $F_* \mathbf{m}' \ne \mathbf{m}'$. However, as the measure $\mathbf{m}_+'$ is $F_+$-invariant, all the iterations $F_*^n \mathbf{m}'$ project to the measure $\mathbf{m}_+'$ on the skew product over the shift on space of one-sided sequences. Hence, projection on $\mcI$ of any of these iterations gives us the measure~$\mu'$ (which is the projection of $\mathbf{m}_+'$).




By construction, the measure $\mathbf{m}'$ projects to the measure $\nu$ on the base. Due to Remark~\ref{r:lim}, the iterations $F_*^n\mathbf{m}'$ of the measure $\mathbf{m}'$ converge to the measure $\mathbf{m}_{\varphi}$. On the other hand, all these iterations project to the measure $\mu$ on $\mcI$, thus the same holds for their limit. Hence, the measure $\mu$ is the projection on $\mcI$ of the measure $\mathbf{m}_{\varphi}$. This projection is the distribution $\varphi_*\nu$ of the values of the map~$\varphi$, that corresponds to the graph part of the maximal attractor. As we supposed only that  the measure $\mu'$ is a stationary measure on the domain $\mcI_{\mu}$, the stationary measure on this domain is unique (and hence it is the measure $\mu$).

\end{proof}


\section{Smooth Stochastic Perturbations and the Baxendale's Theorem}    \label{s:bax}

In this Section, we complete the proof of Theorem~\ref{t:random}. The conclusions~\ref{t:random:1} and~\ref{t:random:2} of the theorem follow from Lemmas~\ref{l:MFixed}, \ref{l:Markov-intervals}, and~\ref{l:Markov-uniqueness}. To prove that the (random) Lyapunov exponent of any ergodic stationary measure is negative, we introduce stochastic skew products over discrete Markov chains and prove for them Theorem~\ref{t:baxendale} which is an analogue of the Baxendale's~\cite[Theorem 4.2]{Baxendale1989}.

\subsection{Stochastic Skew Products}

In this paper, a stochastic skew product over a Markov chain is an object of 3 components:
\begin{enumerate}
  \item a discrete transitive Markov chain with $N$ states $\{1,\dots, N\}$, given by its transition matrix $(\pi_{ij})$ and the stationary distribution $p_i$ of this chain;
  \item a compact manifold~$M$;
  \item $N$ random variables~$\xi_{ij}$ taking values in~$\DiffIm^1(M)$.
\end{enumerate}
Observe that if all~$\xi_{ij}$ are almost surely constant, then the stochastic skew product is a usual deterministic skew product.

Similarly to Section~\ref{s:random_markov}, denote~$\mathcal{M} := \{ 1, \dots, N \} \times M$ and $\mathcal{M}_i := \{ i \} \times M$. A probability measure~$\mu$ on $\mathcal{M}$ is called stationary if
\begin{equation}\label{e:stoch_img_2}
  \sum_i \pi_{ij}  \cdot \bbE (\xi_{ij})_* \mu_i = \mu_j,
\end{equation}
where~$\mu_i = \mu|_{\mathcal{M}_i}$. The expectation is taken over the $\xi_{ij}$-images of the measure $\mu_i$ (not the ``spacewise'' expectation along $M$).

Assume that for any~$i,j$ we have
\begin{equation}\label{e:lognorm}
\bbE \log \norm{D\xi_{ij}} < +\infty.
\end{equation}
Then for any ergodic stationary probability measure~$\mu$ its random volume Lyapunov exponent (which is the sum of its random Lyapunov exponents) is well defined and is given by the following formula:
\begin{equation}\label{e:lyap-vol}
    \lambda_{vol} = \sum_{i,j} \pi_{ij} \cdot \bbE \int_M \log \Jac \xi_{ij}|_x d\mu_i(x).
\end{equation}

\begin{Thm} [Baxendale's Theorem for stochastic skew products]\label{t:baxendale}
Let a stochastic skew product over a Markov chain~$\Pi$ satisfy~\eqref{e:lognorm}. Then at least one of the following options is true:
\begin{enumerate}
  \item\label{i:esm} there exists an ergodic stationary measure $\mu$ with $\lambda_{vol}(\mu) < 0$;
  \item\label{i:asi} there exists an almost surely $\xi_{ij}$-invariant measure~$\nu$ on~$\mathcal{M}$:
  $$
    \forall i,j \quad (\xi_{ij})_* \nu_i = \nu_j \quad \text{with the probability~$1$}.
  $$
\end{enumerate}
\end{Thm}

\begin{proof}

To prove Theorem~\ref{t:baxendale}, we follow the spirit of the argument
from~\cite[Theorem 4.2]{Baxendale1989}. For any two probabilities~$m_1, m_2$ on some Polish space $X$ such that $m_1 \ll m_2$, the \emph{relative entropy} (also known as the Kullback-Leibler information divergence) of $m_1$ with respect to $m_2$ is
\begin{equation}\label{e:rel-entropy}
h(m_1|m_2) := \sup_{\psi \in C(X)} \left[ \log \int e^\psi \,dm_1 - \int \psi \,dm_2 \right] =
\left\{
  \begin{array}{ll}
    \int \left( \log \frac{dm_1}{dm_2} \right) \,dm_1, & \hbox{if $m_1 \ll m_2$;}   \\
    +\infty, & \hbox{otherwise.}
  \end{array}
\right.
\end{equation}
By definition, $h(m_1|m_2) \ge 0$. Moreover, because $\log x$ is a convex function, $h(m_1|m_2) = 0$ iff $m_1 = m_2$.

A key idea of~\cite{Baxendale1989} is the connection between the relative entropy and the volume Lyapunov exponent for smooth invariant measures.
\begin{Prop}
Let~$\mu$ be an ergodic stationary measure absolutely continuous w.r.t. the Lebesgue one with a bounded density. Then
\begin{equation}\label{e:lyap-vol-entropy}
    \lambda_{vol}(\mu) = -\sum \pi_{ij} \cdot \bbE h ((\xi_{ij})_* \mu_i | \mu_j).
\end{equation}
\end{Prop}

\begin{proof}
By definition,
\begin{multline}\label{e:lyap-rel}
h ((\xi_{ij})_* \mu_i | \mu_j) = \int_{M_j} \log \left. \frac{d (\xi_{ij})_* \mu_i}{d\mu_j}\right|_{y}  \, d (\xi_{ij})_* \mu_i(y)= \int_{M_i} \log \left(\frac {\rho_{\mu_j}(\xi_{ij}(x))}{(\Jac \xi_{ij})(x) \cdot \rho_{\mu_i}(x)}\right) \, d \mu_i(x)=
\\
= -  \int_{M_i} \log \Jac \xi_{ij}|_{x}  \, d \mu_i(x) + \left(
\int_{M_i} \log \left({\rho_{\mu_j}(\xi_{ij}(x))}\right) \, d \mu_i(x) -
\int_{M_i} \log \left({\rho_{\mu_i}(x)}\right) \, d \mu_i(x)\right).
\end{multline}
where $\rho_{\mu_i}$ and $\rho_{\mu_j}$ stay for the densities of $\mu_i$ and $\mu_j$ respectively, and the first equality is due to the change of variables $y=\xi_{ij}(x)$. The latter two integrals are well defined because the density $\rho_{\mu}$ is bounded, and the manifold $M$ is compact.

Now we apply $-\sum_{i,j} \pi_{ij}\cdot \bbE (\cdot)$ to the right hand side of~\eqref{e:lyap-rel}. The first summand becomes~\eqref{e:lyap-vol}. The second one becomes zero because $\mu$ is stationary (see~\eqref{e:stoch_img_2}):
\begin{multline}\label{e:lyap-corr}
\sum_{i,j} \pi_{ij} \bbE \int_{M_i} \log \rho_{\mu_j}(\xi_{ij}(x))\, d\mu_i(x)  -  \sum_{i,j} \pi_{ij} \bbE \int_{M_i} \log \rho_{\mu_i}(x) \,d\mu_i(x) = \\
=\sum_{j} \int_{M_j} \log \rho_{\mu_j}(y)\, d\left(\sum_i \pi_{ij} \bbE (\xi_{ij})_{*} \mu_i\right)(y)   -  \sum_{i} \int_{M_i} \log \rho_{\mu_i}(x) \,d\mu_i(x) = 0.
\end{multline}
This proves~\eqref{e:lyap-vol-entropy}. Note that because $\rho_\mu$ is bounded, $\int_{M_i} \log \rho_{\mu_i}(x) \,d\mu_i(x) < +\infty$. Also $\int_{M_i} \log \rho_{\mu_i}(x) \,d\mu_i(x) > -\infty$ because $M$ is compact. Thus both integrals are finite.
\end{proof}

To apply~\eqref{e:lyap-vol-entropy}, let us smooth the process~$\Pi$. Namely, for any $\eps>0$ fix a random variable $\zeta^\eps$, independent of $\xi_{ij}$, taking values in $\DiffIm^1(M)$, such that the following conditions hold:
\begin{enumerate}
\item for any $x\in M$, the image $\zeta^\eps(x)$ is almost surely within the $\eps$-neighborhood of x
\item the random variables $\zeta^\eps(x)$ are absolutely continuous with the density uniformly bounded by some constant $C^{\eps}$.
\item the norms $\|D\zeta^\eps\|, \|(D\zeta^\eps)^{-1}\|$ are bounded almost surely.
\end{enumerate}

In particular, for the case of $M$ being an interval, $\zeta^{\eps}$ can be chosen as a composition of a contraction with factor $1-\eps$ and a random translation by a shift that is chosen from $[-\eps,\eps]$ with a smooth density.


Denote by~$\Pi^\eps$ the stochastic skew product with~$\xi_{ij}^\eps = \xi_{ij} \ast \zeta^\eps$ where~$\ast$ stands for the convolution. For any~$\eps > 0$, there exists at least one ergodic $\Pi^\eps$-stationary measure~$\mu_\varepsilon$ for~$\Pi^{\eps}$. Note that due to the choice of $\Pi^{\eps}$ this measure automatically has a density that does not exceed the same constant~$C^{\eps}$: indeed, this holds for one-step averaging of any initial measure due to the assumption on the laws of $\zeta^{\eps}(\cdot)$, and the stationary measure is equal to its one-step diffusion.

Now, for any such measure the corresponding volume random Lyapunov exponent (of the new process $\Pi^{\eps}$) is nonpositive, and is zero iff the measure is a.s. invariant. Indeed,
$$
\lambda_{vol}(\mu^{\eps}; \Pi^{\eps})=-\sum_{i,j} \pi_{ij} \cdot \bbE h ((\xi_{ij}^{\eps})_* \mu^{\eps}_i | \mu^{\eps}_j) \le 0,
$$
and due to the properties of the relative entropy the right hand side can be zero only if $(\xi_{ij}^{\eps})_* \mu_i = \mu_j$ almost surely. Now we carefully pass to the limit to obtain the same for the original process $\Pi$.

The space of distributions on~$M$ is compact, so one can extract from the family $\{\mu^{\eps}\}$ a weakly converging subsequence:
$\mu^{\eps_n}\to\mu, \quad \eps_n \searrow 0$ as $n\to \infty$. Then either $\lambda_{vol}^{\eps_n} \to 0$ and thus $h((\xi_{ij}^{\eps_n})_* \mu_i | \mu_j) \to 0$ weakly, or there exists a subsequence $\eps'_n \searrow 0$ such that $\lambda_{vol}^{\eps'_n} \to -\alpha < 0$.

In the latter case,
$$
\sum_{i,j} \pi_{ij} \int \log \Jac \xi_{ij}|_x \, d\mu_i(x) = \lim_{n \to \infty} \sum_{i,j} \pi_{ij} \int \log \Jac \xi_{ij}^\eps|_x \, d\mu_i^\eps(x) = -\alpha < 0,
$$
and hence for at least one ergodic component $\tilde{\mu}$ of~$\mu$ (note that $\mu$ is not necessarily ergodic!) we have
$$
\lambda_{vol} (\tilde{\mu})= \sum_{i,j} \pi_{ij} \int \log \Jac \xi_{ij}|_x \, d\mu_i(x) \le -\alpha < 0.
$$
This is the option~\ref{i:esm}) of the theorem.

On the other hand, the relative entropy is upper semicontinuous as the supremum of continuous functionals. Hence, in the former case for any admissible $(i,j)$ one has almost surely
$h((\xi_{ij})_* \mu_i | \mu_j) = 0$, and thus $(\xi_{ij})_* \mu_i = \mu_j$. This is the option~\ref{i:asi}).

Theorem~\ref{t:baxendale} is proven.
\end{proof}

\subsection{The proof of the Theorem~\ref{t:random}}

Recall the Markov process~$\Pi(F)$ defined in Section~\ref{s:random_markov}. Fix any ergodic stationary measure~$\eta$ of $\Pi(F)$. Restrict~$\Pi(F)$ to the trapping region~$\mcI_\eta$. The process $\Pi(F)$ can be viewed as a stochastic skew product over a Markov chain with
\begin{enumerate}
  \item the Markov chain equal to our given Markov chain in the base;
  \item $M = \mcI_\eta$;
  \item $\xi_{ij}$ equal to the $\delta$-measures at~$f_i \in \DiffIm^1(M)$.
\end{enumerate}

Apply the Baxendale's Theorem to~$\Pi(F)$. The option~\ref{i:asi} is impossible by Corollary~\ref{cor:meas_not_coincide}.
Thus the option~\ref{i:esm} is true, and there exists an ergodic stationary measure in~$\mcI_\eta$ with $\lambda_{vol} < 0$. As $\dim I = 1$, $\lambda_{vol}$ is just its random Lyapunov exponent.

But by Lemma~\ref{l:Markov-uniqueness}, $\eta$ is the unique stationary measure in~$\mcI_\eta$. Hence $\Lyap\eta < 0$. The proof of Theorem~\ref{t:random} is complete.

\section{Proof of the skew products theorem}   \label{s:proof_markov_classif}

%

\begin{Prop}    \label{p:limit_step}
Let the map~$\varphi(\omega)$, $\varphi \colon \Sigma \to S^1$, depend only on finitely many symbols in~$\omega$. Suppose its graph~$\Gamma$ drifts up (down). Then
\begin{enumerate}
\item the pointwise limit of its iterates~$\Gamma_n = F^n(\Gamma)$ as $n \to +\infty$ is the graph of a measurable function $\varphi_{+\infty} \colon \Sigma \to S^1$;
\item the function $\varphi_{+\infty}(\om)$ does not depend on ``future'' of~$\om$:
$$
\forall i \in \bbZ, \ \forall \om = \om_{[-\infty, -1]} \om_{[0, +\infty]}, \om' = \om_{[-\infty, -1]} \om_{[0, +\infty]}', \quad \varphi_{+\infty} (\om) = \varphi_{+\infty} (\om');
$$
\item $\varphi_{+\infty}$ is invariant under~$F$:
$$
\varphi_{+\infty} (\sigma \om) = f_{\om} (\varphi_{+\infty} (\om)).
$$
\end{enumerate}
\end{Prop}

\begin{Rem}
The analogous statement holds for $n \to -\infty$. The limit function $\varphi_{-\infty}$ does not depend on the ``past'' of~$\om$.
\end{Rem}

\begin{proof}
Take any~$\om \in \Sigma$. Let~$\varphi_n$ be the function that corresponds to the graph~$\Gamma_n$. Because~$\Gamma$ drifts up (down), the sequence $(\varphi_n(\om))$ is monotone. It is also bounded because of the invariant repeller graph~$J \subset \Sigma \times (S^1 \setminus I)$. Thus $\varphi_{+\infty}$ is well defined. Because $\varphi_{+\infty}$ is a pointwise limit of continuous functions~$\varphi_n$, it must be measurable.

Now let~$\varphi(\om)$ depend only on~$\om_{-k}, \ldots, \om_m$ in~$\om$. The function~$\varphi_n$ has the form
$$
\varphi_n(\om) = f_{n, \sigma^{-n}\om} (\varphi (\sigma^{-n} \om)).
$$
Then for every~$n\ge m+1$ the function~$\varphi_n$ depends only on~$\om_{-k-n}, \ldots, \om_{-1}$. In particular,~$\varphi_n$ does not depend on~$\om_{[0,+\infty]}$. Thus the limit function~$\varphi_{+\infty}=\lim_{n\to+\infty} \varphi_n$ is also independent of~$\om_{[0,+\infty]}$.
%

The invariance of $\varphi_{+\infty}$ immediately follows from the definition of limit.
\end{proof}

For any measurable function~$\varphi$ denote by~$\mathbf{m}_\varphi$ the lift of the base measure~$\nu$ to the graph~$\Gamma$ of~$\varphi$. The measure~$\mathbf{m}_\varphi$ is invariant provided that~$\Gamma$ is invariant.


\begin{Prop}    \label{p:meas_on_limit}
Suppose a function~$\varphi \colon \Sigma \to S^1$ is independent of the future and the graph of~$\varphi$ is $F$-invariant. Then there exists a stationary measure~$\mu$ of the process~$\Pi(F)$ such that
$\pi_* \mathbf{m}_{\varphi} = \mathbf{m}_+ (\mu)$, where measure~$\mathbf{m}_+ (\mu)$ is defined by Eq.~\eqref{e:inv_measure}. 
\end{Prop}

\begin{proof}
Because $\varphi$ is independent of the future, the projection of the measure~$\mathbf{m}_\varphi$ to $\Sigma_+ \times S^1$ has the following simple form. The restriction of the projection to each cylinder~$\mathcal{C}_k$ it is the Cartesian product of~$\nu_k^+$ and some measure~$\mu_k$ on $I_k$, see~\eqref{e:inv_measure}. These measures~$\mu_k$ constitute the desired stationary measure~$\mu$ on $\mcI$.
%
%
\end{proof}

\begin{Rem} \label{r:lyap}
By definition, the fiberwise Lyapunov exponent of $\mathbf{m}_\varphi$ equals to the random Lyapunov exponent of $\mu$.
\end{Rem}

\begin{Rem}
The analogous statements holds for any~$\varphi$ that is independent of the past. Measure~$\mu$ is stationary for the process~$\Pi(F^{-1})$.
\end{Rem}

Now we are ready to prove Theorem~\ref{t:markov_classification}.
By Theorem~\ref{t:random}, the process~$\Pi(F)$ has finitely many ergodic stationary measures. Their supports are contained in disjoint vertically sorted trapping regions. These regions correspond to the trapping strips in $\Sigma \times S^1$. By Proposition~\ref{p:Markov-bony}, the maximal attractors of these strips are CBGs. Thus any strip has a unique invariant measure projecting to $\nu$ in the base. By Proposition~\ref{p:SRB}, this measure is physical. The claims~\ref{t:markov_classification:1}, \ref{t:markov_classification:2} are established. Proposition~\ref{p:thin_conj} implies claim~\ref{t:markov_classification:2.5}.

Because the fiber maps are monotonous, the complement to the disjoint union of the trapping strips is also the disjoint union of finitely many vertically ordered step strips. The $F$-images of these strips are step w.r.t. the symbol~$\omega_{-1}$ of $\om$. These new strips are inverse trapping.

Let us see what happens in any of these strips. According to Propositions~\ref{p:limit_step} and~\ref{p:meas_on_limit}, the limits of the backward iterates of its lower and upper boundaries are some invariant measurable graphs~$\Gamma_L$ and $\Gamma_U$. They support invariant measures, namely, the lifts of $\nu$. These invariant measures in turn correspond to some stationary measures~$\mathbf{m}_L$, $\mathbf{m}_U$ in~$\Pi(F^{-1})$. Suppose $\mathbf{m}_L$ and $\mathbf{m}_U$ do not coincide. Then there is an invariant attracting graph between them which has its own trapping strip. Thus the inverse trapping strip can decomposed into at least two strips. This contradiction proves that $\mathbf{m}_L = \mathbf{m}_U$. The unique stationary measure~$\mathbf{m}_L = \mathbf{m}_U$ corresponds to the unique repeller within the inverse trapping strip. The claim~\ref{t:markov_classification:4} is proven.

The claim~\ref{t:markov_classification:3} follows from Rem.~\ref{r:lyap}.

Finally, the drifting up (down) border of any trapping strip converges a.e. to the corresponding attractor as $n \to +\infty$ and to the corresponding repeller as~$n \to -\infty$. This makes a.e. point between the attractor and the repeller do the same. Thus the basins of attractors and repellers cover the whole phase space~$\Sigma \times S^1$ with the exception of zero measure set in the base. The claim~\ref{t:markov_classification:0.5} is proven. Theorem~\ref{t:markov_classification} is complete.

We conclude this paper with an estimate of the number of attractors and repellers.
For any~$f$, denote by~$\AFix(f)$ the set of all attracting fixed points of~$f$.
\begin{Prop}    \label{p:att_num}
  Let~$\om \in \Sigma$ be any periodic sequence, $\om = (w) = \dots www \dots$. Let~$\tilde D$ be any trapping strip. Then there exist $a_1, a_2 \in \AFix(f_w)$ (perhaps, $a_1 = a_2$) such that $\Amax(\tilde D) \cap I_\om = [a_1, a_2]$.
\end{Prop}

\begin{proof}
  Indeed, because $\tilde D$ is trapping, $f_w(\tilde D_\om) \subset \tilde D_\om$, and $\Amax (\tilde D) \cap I_\om = \bigcap_{n \ge 0} f_w^n (\tilde D_\om)$ which has to be an interval between some attracting fixed points of~$f_w$.
\end{proof}

\begin{Cor}		\label{c:att_num}
  The number of the attracting CBGs is less or equal to~$\sharp\AFix(f_w)$.
\end{Cor}

\begin{Rem}	\label{r:att_num}
The minimum of $\sharp\AFix (f_w)$, taken over all the periodic sequences $(w)$, is \emph{equal} to the number of trapping strips (and thus to the number of the attracting CBGs).
\end{Rem}

\begin{proof}[Sketch of the proof] First, note that for a generic $\om$, the intersection of the fiber $I_{\om}$ with the trapping strip is contracted exponentially by the dynamics. This is implied by the topological contraction and by the negativity of the Lyapunov exponent. But the skew product is a step one, so this implies that the proportion in the set of all periodic words of period $n$ of words $w$ for which the map $f_w$ is contracting in~$\tilde D \cap I_{(w)}$, tends to 1 as $n$ tends to $+\infty$.
Hence, for $n$ sufficiently large for most sequences of period $n$ the fiber maps over a period have a single attracting fixed point in~$\tilde D$.

In particular, such $\om = (w)$ exists.
\end{proof}


%


\begin{Rem}
  The analogous statements hold for the repelling CBGs.
\end{Rem} 

\section{Acknowledgements}

The authors are very grateful to professor Yu.~Ilyashenko for posing the general question about the existence of ``thick attractors'' which investigation led to the results mentioned above, as well as for fruitful discussions.
The second author would like to thank Universit\'e~de~Rennes~1 and Institut~de~Recher\-che Mathematique~de~Rennes for the hospitality during the initial work on this paper.

\bibliographystyle{plain}

\inputencoding{cp1251}

\end{document}